\documentclass[twocolumn]{autart}

\usepackage{graphicx}
\usepackage{amsfonts,amssymb}
\usepackage{amsmath}
\usepackage{epstopdf}
\usepackage{bm}
\usepackage{dsfont}
\usepackage{cite}
\usepackage{enumerate}
 \usepackage{subcaption}

\newcommand{\openbox}{\leavevmode
  \hbox to.77778em{%
  \hfil\vrule
  \vbox to.675em{\hrule width.6em\vfil\hrule}%
  \vrule\hfil}}

\newtheorem{theorem}{Theorem}
\newtheorem{example}{Example}
\newtheorem{lemma}{Lemma}
\newtheorem{remark}{Remark}

\newtheorem{corollary}{Corollary}
\newtheorem{problem}{Problem}

\newcommand{\diag}{\mathsf{diag}}
\newcommand{\tr}{\mathsf{tr}}
\newcommand{\blkdiag}{\mathsf{blkdiag}}

\newcommand{\sym}{\mathsf{sym}}
\newcommand{\rank}{\mathsf{rank}}

\usepackage{xcolor}

\definecolor{mygreen}{rgb}{0,0.5,0}

\usepackage{soul}

\usepackage{algorithm}
\usepackage{algorithmicx}
\usepackage{algpseudocode}
\graphicspath{{./figures/}}
\allowdisplaybreaks[4]

\begin{document}
\begin{frontmatter}
\title{$\mathcal{H}_{2}$ model reduction for diffusively coupled second-order networks by convex-optimization\thanksref{footnoteinfo}}

\thanks[footnoteinfo]{This paper was not presented at any IFAC
conference. This work was
supported by the National Natural Science Foundation of China Under Project 61761136005, National Natural Science Foundation of China Under Project 62003276 and 61773357, and the Fellowship of Zhejiang Province Postdoctoral Science Foundation Under Project ZJ2020001. The first two authors contributed equally to this work.
}

\author[Paestum,West,Baiae]{Lanlin Yu}\ead{yulanlin1992@gmail.com},    
\author[Rome]{Xiaodong Cheng}\ead{xc336@cam.ac.uk},               
\author[Baiae]{Jacquelien M.A. Scherpen}\ead{j.m.a.scherpen@rug.nl},  
\author[P]{Junlin Xiong}\ead{junlin.xiong@gmail.com}

\address[Paestum]{School of Electrical Engineering and Automation, Hefei University of Technology, Hefei 230009, China.}
\address[West]{Westlake Institute for Advanced Study, Westlake University, Hangzhou 310024, China.}  
\address[Rome]{Department of Engineering, University of
Cambridge, Trumpington Street, Cambridge, CB2 1PZ, United Kingdom.}             
\address[Baiae]{Jan C. Willems Center for Systems and Control, Engineering and Technology Institute Groningen, Faculty of Science and Engineering, University of Groningen, Nijenborgh 4, 9747 AG Groningen, The Netherlands.}        
\address[P]{Department of Automation, University of Science and Technology of China, Hefei 230026.}

\begin{keyword}
Second-order networks, diffusive coupling, $\mathcal{H}_{2}$ model reduction, linear matrix inequality, convex-optimization
\end{keyword}                             

\begin{abstract}
This paper provides an $\mathcal{H}_{2}$ optimal scheme for reducing diffusively coupled second-order systems evolving over undirected networks. The aim is to find a reduced-order model that not only approximates the input-output mapping of the original system but also preserves crucial structures, such as the second-order form,  asymptotically stability, and diffusive couplings. To this end, an $\mathcal{H}_{2}$ optimal approach based on a convex relaxation is used to reduce the dimension, yielding a lower order asymptotically stable approximation of the original second-order network system. Then, a novel graph reconstruction approach is employed to convert the obtained model to a reduced system that is interpretable as an undirected diffusively coupled network. Finally,
the effectiveness of the proposed method is illustrated via a large-scale networked mass-spring-damper system.
\end{abstract}

\end{frontmatter}

\section{Introduction}
Second-order network systems with diffusive couplings are found in a variety of applications, such as mass-spring-damper networks \cite{van2013port}, distributed power grids \cite{dorfler2014sparsity} and electrical circuits \cite{schilders2008model,yan2008RLCK}. With the increasing number of interconnected units in a network, the order of its dynamical model  can easily become high-dimensional, which complicates the analysis and synthesis in the network. It motivates the system approximation for a reduced-order network model that captures the main features of the original one \cite{cheng2021review}. Particularly, for the model reduction problem of second-order networks in this paper, we aim for two goals, namely, approximation of the input-output behavior, and preservation of the network structure with diffusive couplings. The latter essentially requires to restore a Laplacian matrix in the obtained reduced-order model. Such a structure is crucial for  describing the information or energy spreading in networks and hence determines the stability of the entire system \cite{cencetti2018pattern}. Furthermore, consensus, a
widespread phenomenon in networked systems, is also realized based on the diffusive couplings \cite{ren2005survey}, and therefore it is useful to preserve the Laplacian structure for realizing the consensus property in the reduced-order model.

Over the past decades, the study of structure preserving model reduction for network systems has drawn profound interest (see \cite{Monshizadeh2014, besselink2016clustering, ishizaki2015clustered, jongsma2016model, cheng2017reduction,cheng2019directed, necoara2020h_2,cheng2017balancedb} and the references therein). Most of these methods can be classified into two families: clustering-based methods \cite{Monshizadeh2014,  besselink2016clustering, ishizaki2015clustered,ishizaki2013model,jongsma2016model,cheng2019directed, cheng2016model,cheng2017reduction} and balanced truncation methods \cite{cheng2017balanceda,cheng2017balancedb}. The balanced truncation method has been extended to solve the structure preserving model reduction problem for first-order network systems \cite{cheng2017balanceda,cheng2017balancedb}, in which a priori approximation error bound is guaranteed. However, it is not clear how balanced truncation can be applied to second-order network systems.
Although this method have extended to the general second-order case \cite{chahlaoui2006balancing2o, reis2008balanced}, there is no guarantee on either an error bound or network structure. Recently, clustering-based model reduction methods
\cite{ishizaki2015clustered,cheng2016model,cheng2017reduction} have been extended to preserve the network structure for the second-order network systems. However, how to select clusters to achieve the minimal approximation error is an open problem.

In this paper, we focus on convex-optimization techniques, which have already shown satisfactory performances for structure-preserving model reduction problems for e.g., bilinear systems \cite{couchman2011model, qi2016time}, negative imaginary systems \cite{yu2017h, yu2019h}, and input-to-state stable nonlinear systems \cite{ibrir2018projection}. However, for network systems, model reduction methods based on convex-optimization are rarely studied. Although a convex-optimization approach in  \cite{cheng2020Weighting} is proposed to reduce first-order Laplacian dynamics by optimally choosing edge weights in a reduced-order network, there is no direct extension of the result towards second-order networks.

In \cite{wyatt2012issues}, an iterative rational Krylov-based method is presented for reducing second-order systems. However, it
does not guarantee a decrease in the $\mathcal{H}_{2}$ error in each iteration. In contrast to \cite{Monshizadeh2014, cheng2017balancedb, besselink2016clustering, ishizaki2015clustered, jongsma2016model, cheng2017reduction,cheng2018MAS, lanlin2019}, we formulate the model reduction of second-order systems in an optimization framework, which is relaxed as a convex optimization problem, and thus can be efficiently tackled. 
Furthermore, unlike the Riemannian optimization-based approach in \cite{sato2017riemannian} that requires an iterative computation of coupled Lyapunov equations, our method just needs to solve once a linear matrix inequality, which may require a lower computational cost. 
Compared to the method in  \cite{cheng2017balancedb}, a new graph reconstruction method is presented which may produce a network topology that is non-complete.

The rest of this paper is organized as follows. The problem setting is introduced in Section~\ref{sec:preliminaries}, and the main results are presented in Section
\ref{sec: main results}, which 
includes
the convex-optimization approach for reducing second-order systems and a novel graph reconstruction scheme.  In Section~\ref{sec:example},  the proposed method is illustrated by an example and compared with the clustering-based method in \cite{ishizaki2015clustered}. Finally, Section~\ref{sec:conclusion} makes some concluding remarks.

\textit{Notation:} 
The symbol $\mathbb{R}$ denotes the set of real numbers. For a given real matrix $A$, $A^{-1}$ and $A^\top$ stand for the inverse and transpose of $A$, $\sym(A)$ indicates $A^\top+A$, and the columns of $A^{\perp}$ form a basis of the null space of
$A$, that is, $AA^{\perp}=0$. The notation $P>0$~$(\geq0)$ means that a matrix $P$ is positive definite (semi-definite). $I_n$ is the identity matrix of size $n$, and $\mathds{1}_{n}$ represents a vector in $\mathbb{R}^n$ of all ones. $e_{i}$ represents the $i$-th column of $I_n$.

\section{Preliminaries \& problem formulation}\label{sec:preliminaries}

Consider an undirected graph $\mathcal{G}$ that consists of a   node set $\mathcal{V}: = \{1, 2, \cdots, n\}$ and an edge set $\mathcal{E}\subseteq\mathcal{V} \times \mathcal{V}$. $\mathcal{G}$ is \textit{weighted} if each edge, an unordered pair of elements in $\mathcal{V}$, is assigned a positive value (weight). Let $\omega_{ij} > 0$ be the weight of edge $(j,i)$, and $\omega_{ij} = 0$ if $(j,i) \notin \mathcal{E}$. An weighted undirected graph $\mathcal{G}$ can be characterized by the so-called \textit{Laplacian matrix} $L \in \mathbb{R}^{n \times n}$ defined as
\begin{equation} \label{defn:Laplacian}
L_{ij} = \left\{ \begin{array}{ll}
\sum_{j=1,j\ne i}^{n} \omega_{ij} & \quad ~i = j,\\
-\omega_{ij}  & \text{otherwise.}
\end{array}
\right.
\end{equation}
The {Laplacian matrix} $L$ of a connected undirected graph has the following properties:
	(i) ${L}^\top = {L}$ and $L \mathds{1} = 0$;
	(ii) $L_{ij} \leq 0$ if $i \ne j$, and $L_{ij} > 0$ otherwise;
	(iii)$L \geq 0$ and has only one zero eigenvalue.	
Conversely, a real square matrix satisfying the above conditions is the Laplacian matrix of a connected undirected graph.

In this paper, the following second-order network system is studied:
\begin{equation}\label{origi-1}
    \bm{\Sigma}:
\left\{ \begin{split}
	\Ddot{x}+D\dot{x}+K x & = F u, \\
	y&= H x,
\end{split}
\right.
\end{equation}
 with $D\in\mathbb{R}^{n \times n}$, $K\in\mathbb{R}^{n \times n}$ positive definite, called the damping and stiffness matrices, respectively. $F\in\mathbb{R}^{n\times p}$ and $H\in\mathbb{R}^{q\times n}$ are the input and output matrices.
The diffusive coupling among the nodes is represented by an undirected weighted graph, and the stiffness matrix is formed as $K=V+L$, with $L$ a Laplacian matrix, and $V$ a diagonal matrix with non-negative diagonal elements representing self-loops. To ensure $K$ to be positive definite, we require at least one diagonal entry of $V$ being strictly positive. Moreover, we assume a \textit{proportional damping}, i.e.,
\begin{equation} \label{damping}
	D=\alpha I_{n} + \beta K,
\end{equation}
with $\alpha$ and $\beta$ positive scalars. Such a damping is also known as Rayleigh damping or classical damping, which has been studied in various applications \cite{ scruggs2009optimal, gondolo2014characterization}. In this paper, the proportional damping assumption is essential for the reconstruction of a reduced second-order network. There are two key properties of the system $\bm{\Sigma}$:
	(1) $\bm{\Sigma}$ is asymptotically stable owing to the positive definiteness of $D$ and $K$ \cite{Semistable}, and
	(2) both $D$ and $K$ are symmetric and diagonally dominant M-matrices.


A variety of physical networks can be modeled in the second-order form \eqref{origi-1}, such as linearized swing equation in power grids \cite{dorfler2014sparsity}, spatially discretized flexible beams \cite{casella2000modelling} and RLCK circuits \cite{yan2008RLCK}. 
\begin{example}\label{example-1}
A mass-spring-damper network is shown in Fig.~\ref{figure-1}, where each node has the same mass and damping, and the nodes are interconnected by springs.
\begin{figure}[t]
  \centering
  \includegraphics[scale=0.6]{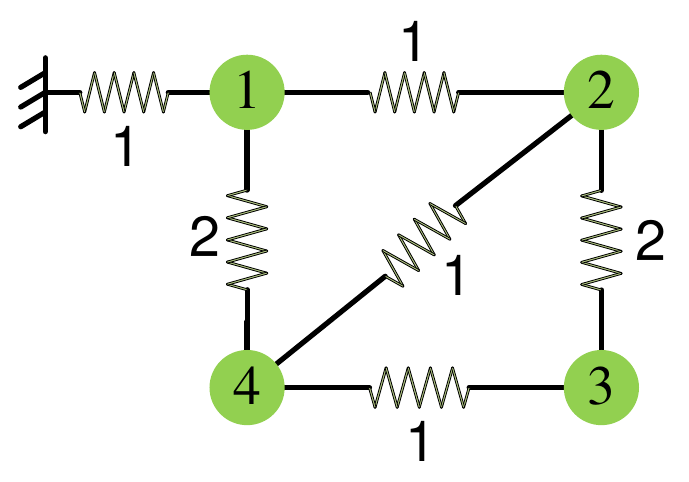}
  \caption{A simple mass-spring-damper network.}
  \label{figure-1}
\end{figure}
The system can be written in the form of \eqref{origi-1} with $D = I$, and
 \begin{equation*}
   K=\begin{bmatrix}
   \begin{smallmatrix}
   4&~ -1&~ 0&~ -2\\
  -1&~  4&~ -2&~ -1\\
   0&~ -2&~ 3&~ -1\\
  -2&~ -1&~ -1&~ 4
  \end{smallmatrix}
   \end{bmatrix}, \quad F=\begin{bmatrix}
   \begin{smallmatrix}
    1\\
    0\\
    0\\
    0
  \end{smallmatrix}
   \end{bmatrix},
 \end{equation*}
 where $K = V+L$ with
 \begin{equation*}
   V=\diag\{1, 0, 0, 0\},~\text{and}~L =\begin{bmatrix}
   \begin{smallmatrix}
    3&~ -1&~  0&~ -2\\
   -1&~  4&~ -2&~ -1\\
    0&~ -2&~ 3&~ -1\\
   -2&~ -1&~ -1&~  4
   \end{smallmatrix}
   \end{bmatrix},
 \end{equation*}
where $L$ is a Laplacian matrix associated with an undirected graph and indicates the strength of diffusive coupling among the nodes connected by the springs in Fig.~\ref{figure-1}.
\end{example}

The model reduction problem for second-order network systems is then formulated as follows.
\begin{problem}\label{oriprobelm}
Given a second-order system $\bm{\Sigma}$ in \eqref{origi-1}, find a reduced second-order network model  
\begin{equation}\label{reduced-1}
    \bm{\Sigma_r}:
\left\{\begin{split}
     	\Ddot{x}_r+D_r\dot{x}_r+K_r x_r & = F_r u, \\
	y_r&= H_r x_r,
\end{split}
\right.
\end{equation}
with $x_r \in\mathbb{R}^{r}$, $y_r\in\mathbb{R}^{q}$, and dimension $1\leq r<n$, such that $\bm{\Sigma_r}$ preserves the network structure, and
the reduction error $\lVert \eta(s)- \eta_r(s) \rVert_{\mathcal{H}_2}$ is as small as possible,  where
$
\eta(s)  = H (s^2 I_n + s D + K)^{-1} F$, and $
\eta_r(s)  = H_r (s^2 I_r + s D_r + K_r)^{-1} F_r.
$
\end{problem}

 We say the reduced-order model is \textit{network structure-preserving}, if $D_r \in\mathbb{R}^{r\times r}$ and $K_r \in\mathbb{R}^{r\times r}$ are positive definite and remain symmetric and diagonally dominant M-matrices. With this structural property, $K_r$ can be written as  $K_r=V_{r}+L_{r}$, where $V_{r}$ is a non-negative diagonal matrix, and $L_{r}$ is an undirected graph Laplacian matrix and thus preserves diffusive couplings among the nodes in the reduced network. This property also applies to the reduced damping matrix $D_r$. In this case, this reduced second-order model preserves the network structure with diffusive couplings.

\section{Main Results}\label{sec: main results}
A two-step approach is presented in this section, where the second-order network system is first reduced by using a convex-optimization approach, and then the resulting reduced-order model is converted into a network system via a graph reconstruction procedure.
\subsection{Model Reduction of Second-Order Systems via Convex Optimization}
\label{sec:H2optimal}

 We make this subsection self-contained. To reduce the interconnected second-order system \eqref{origi-1}, we present an $\mathcal{H}_{2}$ method based on convex optimization. It is worth emphasizing that the method proposed in this subsection is applicable to general second-order systems without the proportional damping assumption in \eqref{damping}.

Before proceeding, the following lemma is provided to characterize the existence of an optimal reduced second-order model of $\bm{\Sigma}_{r}$ in terms of the $\mathcal{H}_{2}$ reduction error.

\begin{lemma}\label{lemma:H_2}
Consider the interconnected second-order system \eqref{origi-1} with positive definite matrices $D$ and $K$.
If there exist positive definite matrices $K_{r}$, $D_{r} \in \mathbb{R}^{r \times r}$, $P \in \mathbb{R}^{2(n+r) \times 2(n+r)}$, and a non-null matrix $F_{r}\in\mathbb{R}^{r\times p}$, such that  the following optimization problem is solvable
\begin{align}\label{opti-1}
 &\min_{P, K_r, D_{r}, F_{r}, H_{r}}~ \tr(HP_{11}
   H^\top-2H_{r}P^\top_{21}
  H^\top+H_{r}P_{31}
   H^\top_{r})
   \nonumber \\
 &\mathrm{s.t.}~~~ P=\begin{bmatrix}
 P_{11}& P_{12}&P_{21}&P_{22}\\
 P^\top_{12}& P_{13}&P_{23}&P_{24}\\
 P^\top_{21}& P^\top_{23}&P_{31}&P_{32}\\
 P^\top_{22}& P^\top_{24}&P^\top_{32}&P_{33}\\
 \end{bmatrix}>0,
 \nonumber \\
 &\qquad PA^\top_{e}+A_{e}P+B_{e}B^\top_{e}=0,
\end{align}
with
\begin{equation}\label{AeBe}
 A_{e}=\begin{bmatrix}
  0&~ I_{n}&~ 0 &~ 0\\
  -K&~ -D&~ 0&~ 0\\
  0&~ 0&~ 0&~ I_{r}\\
  0&~ 0&~ -K_{r}&~ -D_{r}\\
\end{bmatrix}, ~
B_{e}=\begin{bmatrix}
 0\\
 F\\
 0\\
 F_{r}
\end{bmatrix},
\end{equation}
then the output matrix $H_{r}:=HP_{21}P^{-1}_{31}$ minimizes the reduction error $\lVert \eta(s)- \eta_r(s) \rVert_{\mathcal{H}_2}$.
\end{lemma}

\begin{pf}
Let $G_e(s) = C_e(sI - A_e)^{-1} B_e$ with $A_{e}$, $B_{e}$ defined in \eqref{AeBe}, and
$C_{e}=\begin{bmatrix}
H &0 &-H_{r} &0
\end{bmatrix}$. We have $\lVert \eta(s)- \eta_r(s) \rVert_{\mathcal{H}_2} = \lVert G_e(s) \rVert_{\mathcal{H}_2}$. As $D, K>0$, and $D_{r}, K_{r}>0$, the two systems \eqref{origi-1} and \eqref{reduced-1} are asymptotically stable \cite{shieh1987stability}. Therefore, $G_e(s)$ is asymptotically stable, and
\begin{align*}
 \lVert G_e \rVert_{\mathcal{H}_2}^2  &=\tr(C_{e}PC^\top_{e})
 \\
 &= \tr(HP_{11}H^\top
 -2H_{r}P^\top_{21}H^\top+H_{r}P_{31}
 H^\top_{r}),
\end{align*}
According to Propositions 10.7.2 and 10.7.4 \cite{bernstein2009matrix}, the gradient of the above function can be given as
\begin{equation*}
\frac{\partial\tr(C_{e}PC^\top_{e})}{\partial H_{r}}=
-2P^{\top}_{21}H^{\top}+2P_{31}H^{\top}_{r}.
\end{equation*}%
The optimal $H_{r}$ that minimizes $\|G_{e}(s)\|^{2}_{\mathcal{H}_{2}}$ is obtained when the gradient of the above function satisfies $\frac{\partial\tr(C_{e}PC^\top_{e})}{\partial H_{r}}=0$, which
follows that $H_{r}=HP_{21}P^{-1}_{31}$, since $P_{31}>0$.
	\hfill $\openbox$
\end{pf}

\vspace{-20pt}
Lemma~\ref{lemma:H_2} implies that if we can find matrices $D_{r}$, $K_{r}$, $F_{r}$ in \eqref{reduced-1}, and  $H_{r}=HP_{21}P^{-1}_{31}$ satisfying
conditions \eqref{opti-1}, then \eqref{reduced-1} is an optimal reduced-order model in terms of the $\mathcal{H}_2$ norm.
However, finding matrices $D_{r}$, $K_{r}$, $F_{r}$, and $P$ as the optimal solution of the problem \eqref{opti-1} is not straightforward, since the constraints are nonlinear and thus difficult to be tackled numerically. The following theorem is then provided to relax the optimization problem \eqref{opti-1}, which can be solved efficiently as a convex-optimization problem.

\begin{theorem}\label{theorem-reducedasy}
Given the interconnected second-order system \eqref{origi-1}. If there exist matrices $\hat{P}_{11}=\hat{P}^\top_{11}>0$, $\hat{P}_{11}\in\mathbb{R}^{n\times n}$, $\hat{P}_{12}\in\mathbb{R}^{n\times n}$, $\hat{P}_{13}=\hat{P}^\top_{13}>0$, $\hat{P}_{13}\in\mathbb{R}^{n\times n}$, $\hat{P}_{31}=\hat{P}^\top_{31}>0$, $\hat{P}_{31}\in\mathbb{R}^{r\times r}$, a full column rank matrix $\hat{P}_{21}\in\mathbb{R}^{n\times r}$, a scalar $\gamma>0$, such that the following optimization problem is solvable
\begin{subequations}\label{optimization2}
\begin{align}
  \min_{\hat{P}>0}&~~~ \gamma \\
\label{trace1}
  \rm{s.t.}\quad&  \tr\left(H
(\hat{P}_{11}-2X)H^\top\right)<\gamma,   \\
\label{linearinequ1}
&\Pi=\begin{bmatrix}
  \sym(\hat{P}_{12})& \Pi_{12}\\
  \star& \Pi_{22}
\end{bmatrix}<0, \\
\label{linearinequ2}
&\Phi=\begin{bmatrix}
\sym(\hat{P}_{12})& \Phi_{12}\\
\star& \Phi_{22}
\end{bmatrix}<0,\\
\label{linearinequ3}
&\Xi=\hat{P}_{11}-2X>0,\\
\label{hatP1}
 &\hat{P}=\begin{bmatrix}
\hat{P}_{11}& \hat{P}_{12}&\hat{P}_{21}&0\\
\hat{P}^\top_{12}& \hat{P}_{13}&0&0\\
\hat{P}^\top_{21}& 0&\hat{P}_{31}&-\hat{P}_{31}\\
0& 0&-\hat{P}_{31}&  2\hat{P}_{31}\\
\end{bmatrix}>0,
\end{align}
\end{subequations}
where $X=\hat{P}_{21}\hat{P}^{-1}_{31}\hat{P}^\top_{21}$, $\rank(X)\leq r$,
\begin{align*}
\Pi_{12}&=\hat{P}_{13}-\hat{P}_{11}K -\hat{P}_{12}D,\\
\Pi_{22}&=\sym(-K\hat{P}_{12}-D\hat{P}_{13})
 +FF^\top,\\
\Phi_{12}&=-\hat{P}_{11}K-\hat{P}_{12}D+\hat{P}_{13}+2XK,\\
\Phi_{22}&=\sym(-K\hat{P}_{12}-D\hat{P}_{13}),
\end{align*}
then the reduced second-order model
\begin{equation}\label{reduced-hat}
    \bm{\hat{\Sigma}_r}:
\left\{\begin{split}
	\Ddot{\hat{x}}_r+\hat{D}_r\dot{\hat{x}}_r+\hat{K}_r\hat{x}_r & = \hat{F}_r u, \\
	\hat{y}_r&=\hat{H}_r\hat{x}_r,
\end{split}
\right.
\end{equation}
with
\begin{equation}\label{reducedasym2}
\begin{split}
\hat{K}_{r}&=\hat{P}^{-1}_{31}\hat{P}^\top_{21}K
\hat{P}_{21}\hat{P}^{-1}_{31}, ~\quad \hat{F}_{r}=\hat{P}^{-1}_{31}\hat{P}^\top_{21}F, \\
\hat{D}_{r}&= \hat{P}^{-1}_{31}\hat{P}^\top_{21}D
 \hat{P}_{21}\hat{P}^{-1}_{31},~\quad \hat{H}_{r} =H\hat{P}_{21}\hat{P}^{-1}_{31}
\end{split}
\end{equation}
is asymptotically stable. Moreover, the $\mathcal{H}_2$ approximation error has the following upper-bound
\begin{equation}\label{errorbound}
\lVert \bm{\Sigma}- \bm{\hat\Sigma}_{r} \rVert_{\mathcal{H}_2} < \gamma.
\end{equation}
\end{theorem}
The detailed proof is found in Appendix A. Theorem~\ref{theorem-reducedasy} shows that a reduced second-order system \eqref{reduced-hat} can be obtained by solving the optimization problem \eqref{optimization2}, which actually achieves a local optimum that minimizes the $\mathcal{H}_{2}$ reduction error. Compared with the original problem  \eqref{opti-1}, the structure constraint on matrix $\hat{P}$ in the
optimization problem \eqref{optimization2} is more strict, yielding a tighter feasible solution set. Thus, it may not produce an optimal solution to minimize the error $\lVert \bm{\Sigma} -\bm{\hat{\Sigma}_{r}} \rVert_{\mathcal{H}_2}$. Instead, it gives an upper bound $\gamma$ for this error, as given in \eqref{errorbound}.

In Theorem \ref{theorem-reducedasy}, the reduced subspace is captured by $\hat{P}_{21} \hat{P}_{31}^{-1}$, which leads to the reduced second-order model \eqref{reduced-hat} satisfying the following property.
\begin{corollary}
	Consider the interconnected second-order system \eqref{origi-1} with positive definite matrices $D$ and $K$. Then, the reduced second-order model \eqref{reduced-hat} obtained by solving the optimization problem \eqref{optimization2} is asymptotically stable with positive definite matrices $\hat{D}_r$ and $\hat{K}_r$.
\end{corollary}
\vspace{-10pt}
\begin{pf}
	Note that $\hat{P}_{21}$ is imposed to have  full rank, that is $\rank(\hat{P}_{21}) = r$. Thus, the matrix $W: = \hat{P}_{21}\hat{P}^{-1}_{31} \in \mathbb{R}^{n \times r}$ has full column rank with $\rank(W)=\rank(\hat{P}_{21})=r$, due to invertible $\hat{P}^{-1}_{31}$. As a result, the matrices $\hat{D}_r = W^\top D W$ and $\hat{K}_r = W^\top K W$ are positive definite as $D>0$ and $K>0$. The stability of second-order model \eqref{reduced-hat} then follows immediately from \cite{shieh1987stability}.
	\hfill $\openbox$
\end{pf}

The optimization problem \eqref{optimization2} is not convex due to the rank constraints on $\hat{P}_{21}$ and $X$.
Next, we present an numerical algorithm to efficiently solve the optimization problem \eqref{optimization2}, see Algorithm~\ref{alg}.
\begin{algorithm}
\caption{Convex-optimization approach for reducing the interconnected second-order system $\bm{\Sigma}$}
\begin{algorithmic}[1]
\Require
$D$, $K$, $F$, $H$, reduced-order $r$.
\Ensure
      $\hat{D}_{r}$, $\hat{K}_{r}$, $\hat{F}_{r}$, $\hat{H}_{r}$ in \eqref{reducedasym2}.
\State Solve the following convex optimization problem w.r.t.  $\gamma>0$, $\hat{P}_{11}>0$, $\hat{P}_{12}$, $\hat{P}_{13}>0$, and $X_{1}>0$:
\begin{equation}\label{optimizationAlgo}
  \begin{split}
     \min&\quad \gamma \\
     s.t.& \quad X=\blkdiag\{X_{1}, 0\}\geq0, \\
     &\quad \eqref{trace1}-\eqref{hatP1}, ~X_{1}\in\mathbb{R}^{r\times r}.
  \end{split}
\end{equation}
\State Take the Schur decomposition
$X_{1}=UZU^\top$, with a unitary matrix $U$ and quasi-triangular matrix $Z$.
\State Let $\hat{P}_{21}=\begin{bmatrix}
\begin{smallmatrix}
 U\\
 \bm{0}_{(n-r)\times r}
 \end{smallmatrix}
\end{bmatrix}$, $\hat{P}_{31}=Z^{-1}$.
\State Compute $\hat{D}_{r}$, $\hat{K}_{r}$, $\hat{F}_{r}$, $\hat{H}_{r}$ using \eqref{reducedasym2}.
\end{algorithmic}
\label{alg}
\end{algorithm}

Note that the optimization problem \eqref{optimizationAlgo} is convex and thus can be efficiently solved. Moreover, $\hat{P}_{21}$ is guaranteed to have full rank, and $\rank (X) \leq r$. The key ingredient for the algorithm is
a structured $X$ in the form of $\blkdiag\{X_{1}, 0\}$. This consideration is inspired by \cite{ibrir2018projection}, which
deals with linear first-order systems.
With the structured $X$, the equation $X=\hat{P}_{21}\hat{P}^{-1}_{31}\hat{P}^\top_{21}$ is simplified to a Schur decomposition.
Furthermore, $X$ is not unique, as $\hat{P}_{21}$ can be changed as long as $H_{s}\hat{P}_{21}\neq0$ holds.

\begin{remark}\label{remark-2}
Both Theorem~\ref{theorem-reducedasy} and Algorithm~\ref{alg} can be applied to more general second-order systems with a positive definite $K$ and a proportional damping matrix $D$. Moreover, our approach can preserve the proportional damping structure in the reduced-order model, i.e., $\hat{D}_r$ is again a proportional damping matrix. To obtain a better reduced-order model, the Riemannian optimal model reduction method \cite{sato2017riemannian} requires an iterative computation of coupled Lyapunov equations and the optimization of the initial point, which yields a high computational cost if the system dimension is large. Furthermore, the iterative rational Krylov-based method in \cite{wyatt2012issues} does not guarantee a decrease in the $\mathcal{H}_{2}$ error in each iteration. In contrast, our method can obtain a local optimal reduced-order model can be obtained by solving a convex optimization problem.
\end{remark}

\subsection{Reconstruction of diffusive couplings}\label{subseciton-Reconstruction of diffusive}

With Algorithm~1, we obtain the reduced second-order model $\bm{\hat{\Sigma}_{r}}$ as in \eqref{reduced-hat}. However, the matrices $\hat{D}_r$ and $\hat{K}_r$ may not be used to present a network with diffusive couplings, and thus  the reduced-order model as in \eqref{reduced-hat} is not in a network form. In this subsection, we find a reduced-order network model with diffusive couplings that has the same input-output mapping as the reduced second-order system as in \eqref{reduced-hat}.

 Note that the eigenvalues of $\hat{K}_{r}$ are positive real. Thus, $\hat{K}_{r}$ can  be rewritten as
	\begin{equation}\label{K_r}
	\hat{K}_{r}=\lambda_{r}I_{r}+\mathcal{L}_{r},
	\end{equation}
	where $\lambda(\mathcal{L}_{r})=\{\lambda_{1}-\lambda_{r}, \cdots, \lambda_{r-1}-\lambda_{r}, 0\}$ are non-negative real, and $\mathcal{L}_{r}$ has exactly one zero eigenvalue with $\lambda_{1}\geq\lambda_{2}\geq\cdots\geq \lambda_{r-1}>\lambda_{r}>0.$
	However, $\mathcal{L}_{r}$ is not a \textit{Laplacian matrix}, and thus it cannot interpret diffusive couplings. According to \cite[Them. 12]{cheng2017balancedb}, since the eigenvalues of $\mathcal{L}_{r}$ are non-negative real and $\mathcal{L}_{r}$ has exactly one zero eigenvalue, there always exists a {Laplacian matrix} $L_r$ similar to the matrix $\mathcal{L}_{r}$ in \eqref{K_r}, namely, $L_r$ and $\mathcal{L}_{r}$ have the same eigenvalues. This implies that there always exists a linear transformation $K_{r}=U_{r}\hat{K}_{r}U^{\top}_{r}$ such that $K_{r}$ is a stiffness matrix, which represents the diffusive couplings of the reduced second-order system. However, in terms of network reconstruction, \cite{cheng2017balancedb} only provides a procedure to construct a non-sparse graph representation where the vertices in the reduced network are fully connected.

In contrast, this paper provides an alternative graph reconstruction method that may induce a non-complete reduced network. This essentially requires a similarity transformation of $\hat{K}_r$, which results in a matrix $K_{r}$ with the same eigenvalues of $\hat{K}_{r}$ but having a network interpretation. The feasibility of this novel graph reconstruction method is guaranteed in the following theorem.

\begin{theorem}\label{theorem-2}
	Consider any positive definite matrix $\hat{K}_r$ whose eigenvalue decomposition is given as $\hat{K}_r=\mathcal{U}\Lambda\mathcal{U}^{\top}$, with $\hat{\Lambda}=\diag\{\lambda_{1}, \cdots, \lambda_{r-1}, \lambda_{r}\}$. Define a matrix
	\begin{equation}\label{mathcal_V}
	\mathcal{V}=\begin{bmatrix}
	\frac{1}{\sqrt{r}}&~ -\mathds{1}^{\top}_{r-1}T\\
	\frac{1}{\sqrt{r}}\mathds{1}_{r-1}&~ T
	\end{bmatrix},
	\end{equation}
	with
	$T\in \mathbb{R}^{(r-1)\times(r-1)}$   a non-singular matrix satisfying
	\begin{equation}\label{tmatrix}
	TT^{\top}=I_{r-1}-\frac{1}{r} \mathds{1}_{r-1}\mathds{1}^{\top}_{r-1}.
	\end{equation}
	The elements of $T$ fulfill $T_{ij}T_{sj}\leq0$ for $i\neq s$, $j\in\{1, \cdots, r-m-1\}$, and $T_{ij}T_{sj}\geq0$ for $i\neq s$, $j\in\{r-m, \cdots, r-1\}$ with $1\leq m\leq(r-2)$.
	Then, 	\begin{equation}\label{Ur}
	U_{r}=\mathcal{V}\mathcal{U}^{\top}
	\end{equation}
	is a unitary matrix, and
	 $K_r =U_r
	\hat{K}_{r} U^{\top}_r$
 	is a symmetric and diagonally dominant M-matrix.
\end{theorem}
\vspace{-10pt}
\begin{pf}
 We first prove that $U_{r}$ is unitary if matrices $\mathcal{V}$, $\mathcal{U}$ and $T$ are constructed as in Theorem \ref{theorem-2}. It is verified from \eqref{mathcal_V} and \eqref{tmatrix} that $\mathcal{V}\mathcal{V}^{\top}=I_{r}$. Moreover, $\mathcal{U}$ is unitary due to the eigenvalue decomposition of a symmetric matrix $\hat{K}_r$. Therefore, we obtain $U_{r}U^{\top}_{r}=I_{r}$.
Next, we show that $K_r =U_r
\hat{K}_{r} U^{\top}_r$
is a symmetric and diagonally dominant M-matrix.

It follows from \eqref{K_r} and $U_{r}U^{\top}_{r}=I_{r}$ that
\begin{equation}\label{reK_r}
  K_{r}=U_r
 \hat{K}_{r} U^{\top}_r=\lambda_{r}I_{r}+\mathcal{V}\Lambda\mathcal{V}^{\top},
\end{equation}
 which is a symmetric and diagonally dominant M-matrix if the positive semi-definite matrix $\hat{\mathcal{L}} : = \mathcal{V}\Lambda\mathcal{V}^{\top}$ is an undirected graph Laplacian. Note that $\hat{\mathcal{L}}$ shares the same spectrum as
 $\Lambda$,
 and it follows from \eqref{mathcal_V} that
\begin{equation}\label{Lronevector}
    \hat{\mathcal{L}}\mathds{1}_{r}=\begin{bmatrix}
    \begin{smallmatrix}
     \mathds{1}^{\mathsf{T}}_{r-1}T\Lambda_{r-1} T^{\top}
     \mathds{1}_{r-1}&~ \star\\
     -T\Lambda_{r-1}T^{\top}\mathds{1}_{r-1}& ~T\Lambda_{r-1}T^{\top}
     \end{smallmatrix}
    \end{bmatrix}
     \mathds{1}_{r}=0,
\end{equation}
where $\Lambda=\diag\{0, \Lambda_{r-1}\}$ with
$\Lambda_{r-1}=\diag\{\lambda_{1}-\lambda_{r}, \cdots, \lambda_{r-1}-\lambda_{r}\}.$ That means the row and column sums of $\hat{\mathcal{L}}$ are zero.

To further show that $\hat{\mathcal{L}}$ represents an undirected graph Laplacian matrix, then we show that $\hat{\mathcal{L}}$ (i) has all positive diagonal elements and (ii) non-positive off-diagonal entries. The first point is not hard to see, as $T\Lambda_{r-1}T^{\top}$ in \eqref{Lronevector} is strictly positive definite. Now, we prove that the off-diagonal entries of $\hat{\mathcal{L}}$ are either negative or zero.

From the property of $T$ in \eqref{tmatrix}, we obtain that
\begin{equation}\label{matrixT}
 \begin{split}
    &\sum^{r-1}_{j=1} T^{2}_{ij}=1-\frac{1}{r}, \qquad \sum^{r-1}_{j=1, i\neq s} T_{ij}T_{sj}=-\frac{1}{r},\\
    &\sum^{r-1}_{j=1}T^{2}_{ij}+\sum^{r-1}_{s=1, s\neq i}\left(\sum^{r-1}_{j=1}T_{ij}T_{sj}\right)=\frac{1}{r},
 \end{split}
\end{equation}
for any $i, s\in\{1, \cdots, r-1\}$. This further implies that
\begin{equation*}
  \begin{split}
     \sum^{r-m-1}_{j=1, i\neq s}T_{ij}T_{sj}&=-\frac{1}{r}-\sum^{r-1}_{j=r-m, i\neq s}T_{ij}T_{sj}\leq0,  \\
     \sum^{r-1}_{j=r-m, i\neq s}T_{ij}T_{sj}&=-\frac{1}{r}-\sum^{r-m-1}_{j=1, i\neq s}T_{ij}T_{sj}\geq0,
  \end{split}
\end{equation*}
from which, we have
\begin{equation*}\label{lrij}
  \begin{split}
    \sum^{r-1}_{j=1}(\lambda_{j}-\lambda_{r})T^{2}_{ij} & \geq(\lambda_{r-1}-\lambda_{r})
  (1-\frac{1}{r})>0,
  \\
    \sum^{r-1}_{j=1, i\neq s}(\lambda_{j}-\lambda_{r})T_{ij}T_{sj}
    & =\sum^{r-m-1}_{j=1, i\neq s}(\lambda_{j}-\lambda_{r})T_{ij}T_{sj} 
    \\
    & +\sum^{r-1}_{j=r-m, i\neq s}(\lambda_{j}-\lambda_{r})T_{ij}T_{sj}\\
    &\leq -\frac{1}{r}(\lambda_{r-m}-\lambda_{r})<0.
  \end{split}
\end{equation*}
Therefore,  $\hat{\mathcal{L}}_{(i+1)(s+1)}\leq0$ for $i, s\in\{1, \cdots, r-1\}$, $i\neq s$. Moreover, according to \eqref{Lronevector}, it holds that
\begin{multline*}
   \hat{\mathcal{L}}_{(i+1)1}=-\sum^{r-1}_{j=1}(\lambda_{j}-\lambda_{r})T^{2}_{ij}  \\
    -\sum^{r-1}_{s=1, s\neq i}\left(\sum^{r-1}_{j=1, i\neq s} (\lambda_{j}-\lambda_{r})T_{ij}T_{sj}\right),
\end{multline*}
which leads to
\begin{equation*}
 \begin{split}
    \hat{\mathcal{L}}_{(i+1)1}&\leq-(\lambda_{r-1}-\lambda_{r})(1-\frac{1}{r})+\frac{(r-2)}{r}(\lambda_{r-1}-\lambda_{r})\\
    &=-\frac{1}{r}(\lambda_{r-1}-\lambda_{r})<0.
 \end{split}
\end{equation*}
As a result, we have shown that $\hat{\mathcal{L}}_{ii}>0$, $\hat{\mathcal{L}}_{ij}\leq0,~\forall~i \ne j$, implying that $\hat{\mathcal{L}}$ in \eqref{reK_r} is regarded as a Laplacian matrix associated with an undirected weighted graph.
This further yields $K_{r}$ as a symmetric and diagonally dominant M-matrix.
\hfill $\openbox$
\end{pf}

With the matrix $U_r$, the transformed matrix $K_{r}$ possesses the structural property that allows $K_{r}$ to be interpreted as an undirected weighted network with  the diffusive couplings. In this sense, a reduced graph can be reconstructed. Besides, there is a freedom in constructing $U_{r}$ by choosing different $T$. By this means, a sparse $K_{r}$ may be obtained with a particular $T$ under some constraints, see Example \ref{example-2}. 
But we should note that it does not always find non-complete graphs with this approach. Whether we can succeed to find a non-complete graph or not is determined by the prescribed eigenvalues.

Theorem~\ref{theorem-2} shows a sufficient condition for  $T$ to produce a Laplacian matrix, but it does not explicitly state how to choose  $T$, particularly to have zeros in the new stiffness matrix $K_r$. We suggest an ad hoc algorithm to do so. Suppose that we intend to enforce $K_r^{(ij)} = K_r^{(ji)} = 0$. Then, a nonlinear constraint is formed as
\begin{align} \label{eq:elemcons}
       e_i^\top K_r e_j=
    e_i^\top (\lambda_{r}I_{r}+\mathcal{V}\Lambda\mathcal{V}^{\top}) e_j = 0,
\end{align}
where $e_i$ denotes the $i$-th column of the identity matrix. Then a set of nonlinear equations is obtained by combing \eqref{eq:elemcons} and \eqref{tmatrix} in Theorem~\ref{theorem-2}. Note that this set of equations does not always give a solution, depending on the prescribed eigenvalues and how many zero elements are enforced. But when it is solvable, we obtain a non-complete reduced graph, as illustrated in Example~\ref{example-2}.

\begin{remark}\label{remark-3}
Note that we may also use the Householder transformation to construct a tridiagonal $K_{r}$. It has been shown in \cite{ishizaki2013model} that there exists a unique Householder transformation $U_r$ such that $K_{r} = U_r \hat{K}_{r} U_r^\top$ becomes a symmetric tridiagonal M-matrix. However, this tridiagonal $K_{r}$ is not necessary diagonally dominant. Although we can write $K_{r} = V_{r}+L_r$ with $L_r$ representing an undirected chain graph, the diagonal matrix $V_{r}$ may contain negative elements, which losses a physical interpretation.
\end{remark}

In the following example, we demonstrate how to implement our graph reconstruction method in Theorem~\ref{theorem-2}, which is compared with the one in \cite{cheng2017balancedb} and the Householder transformation in \cite{ishizaki2013model}.

\begin{example}\label{example-2}
Let $\{0, 0.4384, 2, 4.5616, 7\}$ be the prescribed eigenvalues, and we aim to create a diagonally dominant M-matrix $K_{r}$ whose eigenvalues match the prescribed ones, and $K_{r}$ has some zero elements, indicating a non-complete graph.
Suppose $K^{(14)}_r = K^{(34)}_r=0$. By solving \eqref{tmatrix} in Theorem~\ref{theorem-2} and the following equations
\begin{align*}
    e_1^\top (\lambda_{r}I_{r}+\mathcal{V}\Lambda\mathcal{V}^{\top}) e_4 = e_3^\top (\lambda_{r}I_{r}+\mathcal{V}\Lambda\mathcal{V}^{\top}) e_4=0
\end{align*}
with $\lambda_{r} = 0.4384$, we obtain a solution as
\begin{equation*}
T=\begin{bmatrix}
\begin{smallmatrix}
-0.8235&~  0&~   0.2682\\
0.1481&~   0.7071&~   -0.4776\\
0.5273&~  0&~  0.6870
\end{smallmatrix}
\end{bmatrix},
\end{equation*}
which leads to 
\begin{equation} \label{eq:KrEx}
K_r =\begin{bmatrix}
\begin{smallmatrix}
3&   -1&   -1.5614&         0\\
-1&    5&   -1&   -2.5615\\
-1.5614&   -1&    3&         0\\
0&   -2.5615&         0&    3
\end{smallmatrix}
\end{bmatrix}.
\end{equation}
Moreover, $K_r$ represents an undirected network with diffusive couplings and the topology is shown in Fig. \ref{figure-2}.

\begin{figure}[t]
\begin{minipage}[t]{0.5\linewidth}
	\centering
	\includegraphics[width=0.55\textwidth]{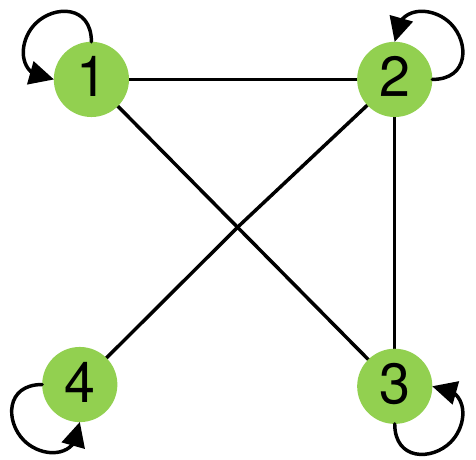}
	\subcaption{}
	\label{figure-2}
\end{minipage}%
\begin{minipage}[t]{0.5\linewidth}
	\centering
	\includegraphics[width=0.55\textwidth]{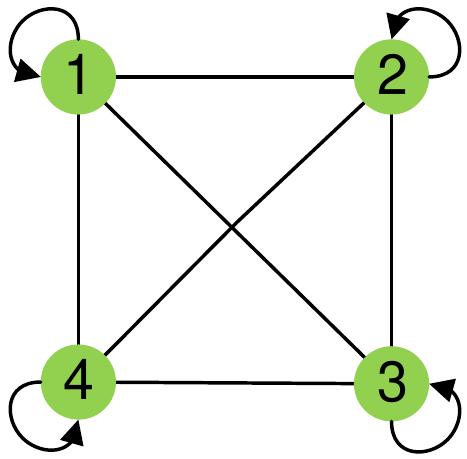}
	\subcaption{}
	\label{figure-21}
\end{minipage}
\caption{(a) The undirected graph corresponding to $K_{r}$ in \eqref{eq:KrEx}. (b) The complete graph   corresponding to $\bar{K}_{r}$ obtained by the procedure in \cite{cheng2017balancedb}.}
\end{figure}

For comparison, we implement the procedure used in the proof of \cite[Them. 12]{cheng2017balancedb} and obtain an alternative graph representation of $\hat{K}_r$ as $\bar{K}_r=0.4384 I_{r}+ \bar{L}_r$ with $\bar{L}_r$ a Laplacian matrix
\begin{equation*}
\bar{L}_r=\begin{bmatrix}
\begin{smallmatrix}
4.5365&~   -1.2443&~  -0.3904&~   -2.4635\\
-1.2443&~    3.3173&~   -0.3904&~   -1.2443\\
-0.3904&~   -0.3904&~    1.6096&~   -0.3904\\
-2.4635&~   -1.2443&~   -0.3904&~    4.5365
\end{smallmatrix}
\end{bmatrix}
\end{equation*}
that represents a complete graph, see Fig. \ref{figure-21}. Furthermore, we use the Householder transformation suggested in \cite[Them. 1]{ishizaki2013model} and obtain a tridiagonal matrix
\begin{equation*}
   \tilde{K}_r=\begin{bmatrix}
   \begin{smallmatrix}
   5&~  -2.4495&~   0&~ 0\\
   -2.4495&~    2.6667&~   -1.8856&~  0\\
   0&~  -1.8856&~    4.3333&~ 0\\
   0&~  0&~    0&~   2
   \end{smallmatrix}
    \end{bmatrix},
\end{equation*}
which is not diagonally dominant and thus losses a network interpretation.
\end{example}


	For the reduced-order model $\bm{\hat{\Sigma}_{r}}$ in \eqref{reduced-hat}, which has the proportional damping, i.e. $\hat{D}_{r}=\alpha I_{r}+\beta\hat{K}_{r}$ with $\alpha>0$, $\beta>0$. Applying the coordinate transformation $x_r=U_r\hat{x}_{r}$ to
$\bm{\hat{\Sigma}_r}$ then leads to a reduced second-order model $\bm{\Sigma_r}$ in the form of \eqref{reduced-1} with
\begin{align}\label{reduced-network}
K_r&=U_r
\hat{K}_{r} U^{\top}_r,~  D_r=U_r
\hat{D}_{r} U^{\top}_r,
\nonumber \\
F_r&=U_r\hat{F}_{r},\qquad H_r=\hat{H}_{r}U^{\top}_r.
\end{align}
Recall the proportional damping assumption in \eqref{damping}, the obtained ${D}_{r}$ from the transformation will be  ${D}_{r}=\alpha I_{r}+\beta {K}_{r}$.
Thus, the reduced model with coefficient matrices in \eqref{reduced-network} possesses the same structure as the original second-order network $\bm{\Sigma}$ in \eqref{origi-1}, and it can be interpreted as a second-order network with reduced number of nodes. Furthermore, the approximation error between the systems $\bm{\Sigma}$ and $\bm{\Sigma_r}$ is evaluated as follows.
\begin{theorem}\label{theo-2}
Consider the original diffusively coupled second-order network $\bm{\Sigma}$ in \eqref{origi-1} and its reduced second-order network model $\bm{\Sigma_r}$ with the matrices in \eqref{reduced-network}. Then, 
we have $\lVert \bm{\Sigma}-\bm{\Sigma_r} \rVert_{\mathcal{H}_2}<\gamma$,
where $\gamma$ is the scalar in \eqref{trace1}.
\end{theorem}

This result follows immediately from that $\bm{\Sigma_r}$ with the matrices in \eqref{reduced-network} is obtained by the coordinate transformation from $\bm{\hat{\Sigma}_r}$ in \eqref{reduced-hat}, and thus they have the same input-output transfer matrices, and $\bm{\Sigma_r}$ is also a solution of Problem~\ref{oriprobelm}.

Although this paper focuses on asymptotically stable second-order network systems, the proposed method can also be easily extended to semi-stable networks studied in \cite{cheng2016model,cheng2017reduction,lanlin2019ecc}, where $K$ is positive semidefinite. The extension can be made by using a system separation as in \cite{besselink2016clustering,cheng2017balancedb}. Taking into account the kernel space of $K$, we have the following decomposition $K = S~ \blkdiag \{ \mathbf{0}_m, \bar{K} \}~ S^\top$, where
 $S$ is unitary, and $m$ is the algebraic multiplicity of the zero eigenvalues of $K$. Here, $S$ can be partitioned as $S = \begin{bmatrix}
 S_0 & S_1
 \end{bmatrix}$
 with 
$K S_0 = 0$.
By defining $z = S^{-1} x = [
z_a \ \
z_s^{\top}
]^{\top}$, with $z_a\in\mathbb{R}^m$ and $z_s\in\mathbb{R}^{n-m}$, the original system \eqref{origi-1} is decomposed into two parts:
\begin{equation}\label{aver-1}
			\ddot{z}_{a}+\alpha\dot{z}_{a}=S_0^{\top}  F u,
			\quad
			y_{a}=  H S_0 z_{a},
\end{equation}
and
\begin{equation}\label{asym-t1}
			\ddot{z}_{s}+\bar{D}
			\dot{z}_{s}+\bar{K}z_{s}=S_1^\top F u,
			\quad
			y_{s}=HS_1 z_{s},
\end{equation}
where $\bar{D} = S_1^\top D S_1$, and $\bar{K} = S_1^\top K S_1$ are positive definite, implying that the system \eqref{asym-t1} is asymptotically stable. By using the proposed $\mathcal{H}_{2}$ optimal model reduction approach in Section~\ref{sec:H2optimal}, we can obtain a reduced second-order model for the system \eqref{asym-t1}. Then combining this reduced  model with the system \eqref{aver-1} results in semi-stable reduced model in the second-order form.
Note that the proportional damping is retained in \eqref{asym-t1}, i.e., $\bar{D} = \alpha I + \beta \bar{K}$. Thereby, the graph reconstitution in Theorem~\ref{theorem-2} can be applied to restore a interconnection structure of diffusive couplings in the reduced  model.

\section{Illustrative Example}\label{sec:example}

In this section, we demonstrate the effectiveness of the proposed model reduction method through an example of complex networks.

\begin{figure}[t]
	\centering
	\includegraphics[width=6cm]{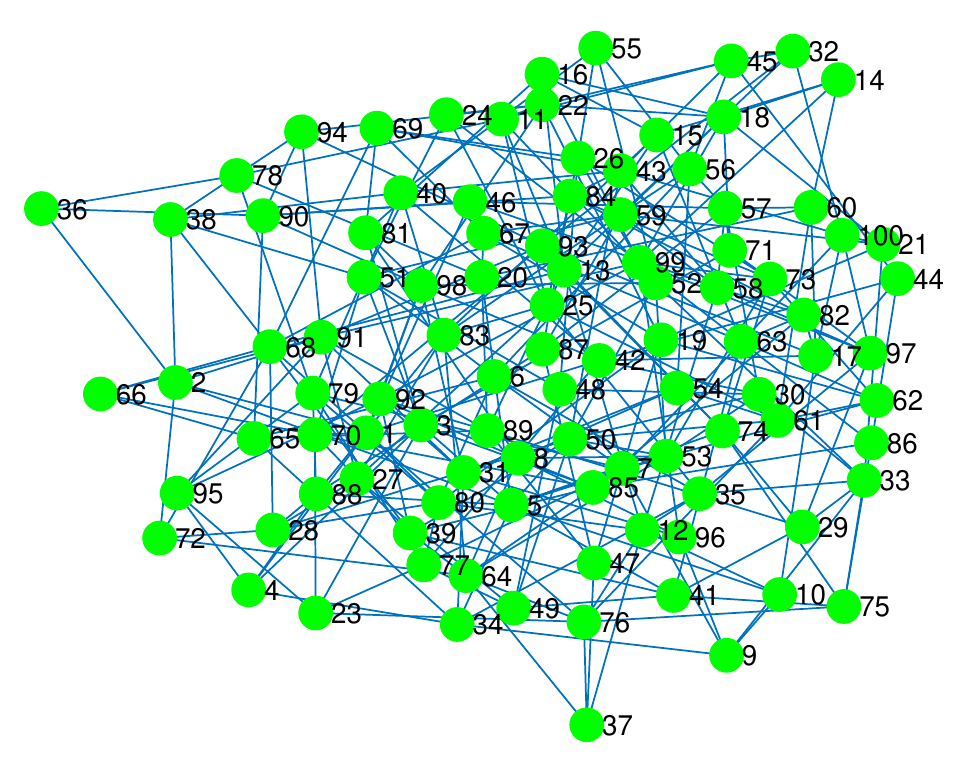}
	\caption{Interconnection topology of the original second-order network (100 nodes).}
	\label{figure-3}
\end{figure}

For comparison, we borrowed the following second-order network in \eqref{origi-1} evolving over the Holme-Kim model composed of 100 nodes \cite{ishizaki2015clustered}, and the interconnection topology is shown in Fig. \ref{figure-3}. In this paper, we select the stiffness matrix $K\in \mathbb{R}^{100\times 100}$ as
\begin{equation*}
K_{ij}=\begin{cases}
 1-\sum^{100}_{j=2} K_{1,j}, & \quad i=1; \\
 -\sum^{100}_{j=1, j\neq i} K_{i,j}, & \quad i\neq 1.
\end{cases}
\end{equation*}
and a proportional damping as $D=\alpha I_{100}+\beta K$ with $\alpha=0.97$ and $\beta=0.15$.
The output and output matrices are chosen as  $F=e^{100}_{1}$ and $H=K-\diag\{e^{100}_{1}\},$ respectively.

\begin{figure}[t]
	\centering
	\includegraphics[width=6cm]{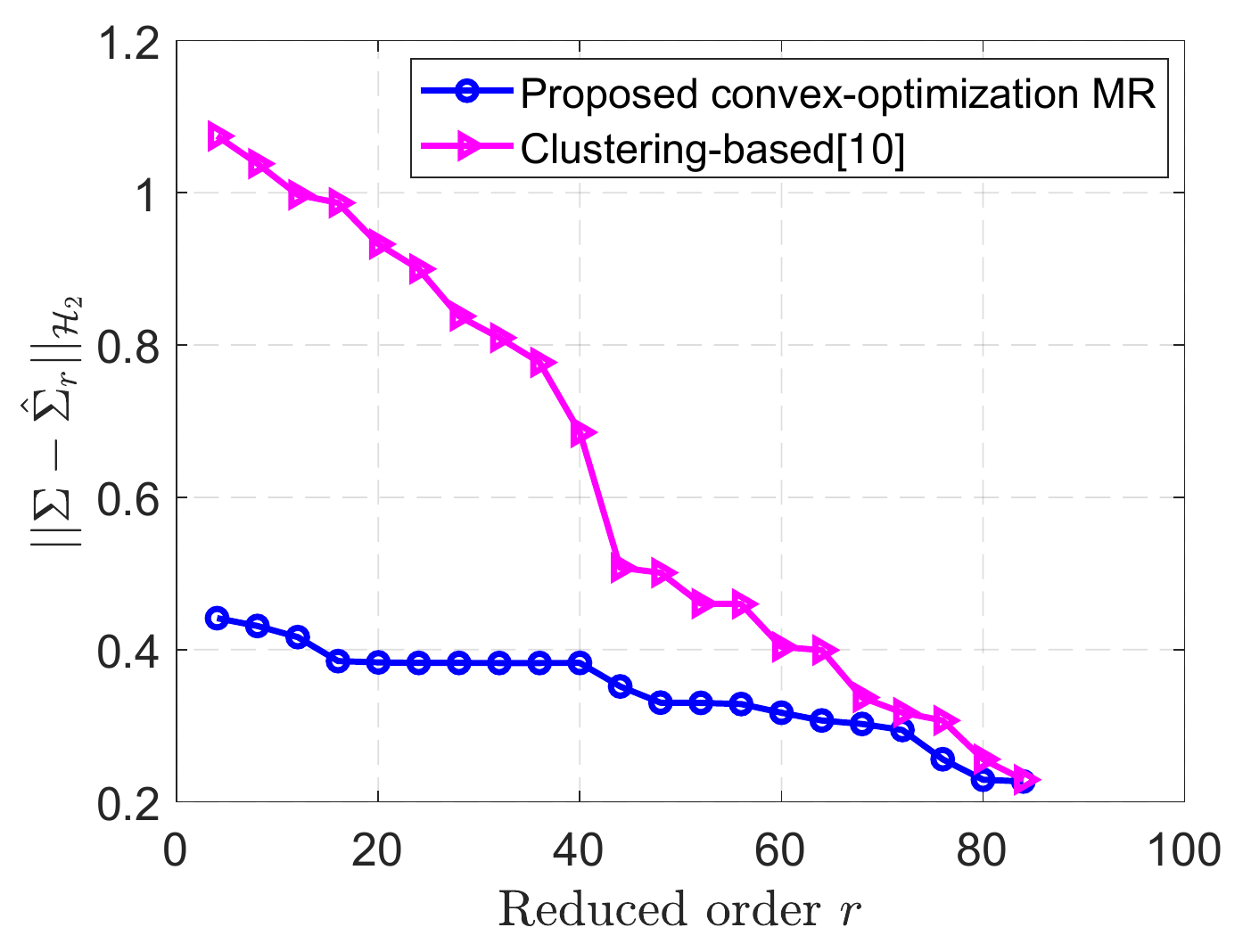}
	\caption{ $\mathcal{H}_{2}$ approximation errors obtained by the proposed model reduction method and clustering-based model reduction method \cite{ishizaki2015clustered}.}
	\label{figure-error-1}
\end{figure}

We reduce the dimension of the second-order network system by two different methods, the clustering-based model reduction method in \cite{ishizaki2015clustered} and the proposed convex-optimization based model reduction method in this paper. Moreover, the reduced-order ranges from 4 to 84 with increments of 4. The $\mathcal{H}_{2}$-norm of the original network system is $1.2661$, and the $\mathcal{H}_{2}$ approximation errors between the original system and the reduced second-order models obtained by Algorithm 1 and the  method \cite{ishizaki2015clustered} are shown in Fig. \ref{figure-error-1}. It can be seen from Fig. \ref{figure-error-1} that the obtained reduced second-order model can approximates the original second-order system well and the $\mathcal{H}_{2}$ approximation error decay as the order of the reduced second-order model increases. Moreover, the proposed  method preserves the second-order network structure and achieves smaller approximation error.

\begin{figure}[t]
\begin{minipage}[t]{0.5\linewidth}
	\centering
	\includegraphics[width=0.58\textwidth]{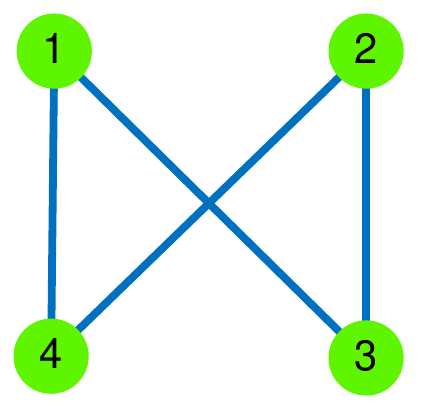}
	\subcaption{}
	\label{figure-5}
\end{minipage}%
\begin{minipage}[t]{0.5\linewidth}
	\centering
	\includegraphics[width=0.58\textwidth]{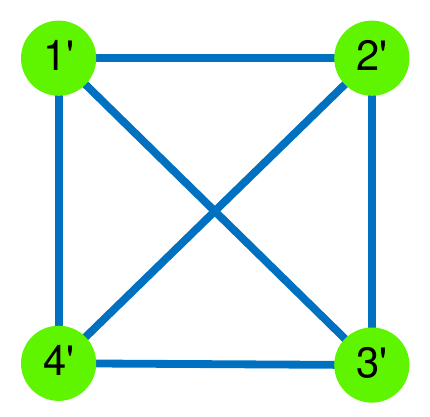}
	\subcaption{}
	\label{figure-6}
\end{minipage}
\caption{(a) The undirected graph corresponding to $K_{4}$ in \eqref{sparse-K4}. (b) The complete graph corresponding to $\tilde{K}_{4}$ in \eqref{nsparse-K4}.}
\end{figure}

To illustrate the effectiveness of our network reconstruction procedure, we consider the obtained reduced model with dimension $4$ as an example, which has the $\mathcal{H}_{2}$ approximation error equal to $0.4371$, and
the eigenvalues of $\hat{K}_{4}$ are given by
$
\lambda(\hat{K}_{4}) =\{9.5631, 7.727, 5.1027, 4.1776\}.$
Based on Theorem \ref{theorem-2}, we select $$T=\begin{bmatrix}
\begin{smallmatrix}
 0.5&~  -0.1845&~  0.6826\\
 0.5&~   0.1845&~  -0.6826\\
-0.5&~  -0.6826&~  -0.1845
 \end{smallmatrix}
\end{bmatrix},$$
which leads to a sparse Laplacian matrix as
\begin{equation}\label{sparse-K4}
K_{4}=\begin{bmatrix}
\begin{smallmatrix}
 6.246&~   0&~  -0.4625&~   -1.6059\\
 0&~    7.039&~   -2.3989&~   -0.4625\\
 -0.4625&~   -2.3989&~    7.039&~   0\\
 -1.6059&~   -0.4625&~   0&~   6.2460
\end{smallmatrix}
\end{bmatrix},
\end{equation}
that has the same spectrum as $\hat{K}_{4}$, and the corresponding interconnection topology is shown in Fig. \ref{figure-5}. Alternatively, we can choose a different $T$ matrix as
$$T=\begin{bmatrix}
\begin{smallmatrix}
0&~  0.309&~   -0.809\\
-0.809&~  0&~   0.309\\
 0.309&~  -0.809&~   0
\end{smallmatrix}
\end{bmatrix},$$
which then yields
\begin{equation}\label{nsparse-K4}
\bar{K}_{4}=\begin{bmatrix}
\begin{smallmatrix}
6.6426&~   -1.6301&~    0.4579&~   -1.2928\\
-1.6301&~    8.0412&~   -1.3463&~   -0.8873\\
0.4579&~   -1.3463&~    5.2973&~   -0.2313\\
-1.2928&~   -0.8873&~   -0.2313&~    6.5889
\end{smallmatrix}
\end{bmatrix},
\end{equation}
representing a complete network with  interconnection topology  shown in Fig. \ref{figure-6}.

%



It can be concluded that the reduced second-order model obtained by the proposed convex-based optimization approach can approximate the original network well. Moreover, a sparse $K_{r}$ may be obtained by using the similarity transformation proposed in Theorem \ref{theorem-2}. That is, a Laplacian matrix associated to an non-complete graph with sparse interconnection is obtained.

Moreover, it can be verified that $K_{4}$ in \eqref{sparse-K4} can be rewritten as
$K_{r}=U_{1}\diag\{9.5631, 7.727, 5.1027, 4.1776\}U^{\top}_{1}$
with a unitary matrix $U_{1}$, which implies  $K_{r}=U_{1}U_{2}\hat{K}_{r}U^{\top}_{2}U^{\top}_{1}$
with $U_{r}=U_{1}U_{2}$, $U_{r}U^{\top}_{r}=I$. Thus, by applying the coordinate transformation $\hat{x}_{r}=U_{r}x_{r}$ to the obtained 4-order model, a second-order network system with diffusive couplings can be obtained, and the interconnection topology is shown in Fig. \ref{figure-5}.

\section{Conclusion}\label{sec:conclusion}

We have developed a novel convex-optimization-based $\mathcal{H}_{2}$ model reduction method for diffusively coupled second-order network systems. A numerical algorithm has been developed to find a local
optimal reduced second-order model. It is worth emphasizing that this algorithm is computationally efficient, as it is constrained by only linear matrix inequalities that can be directly solved  by using efficient convex optimization toolboxes. In addition, by using a new similarity transformation that provided in this paper, the resulting reduced second-order model can be interpreted as an undirected network with diffusive couplings. The main advantage of the proposed method is that a local optimal reduced-order system can be guaranteed in the sense of minimizing the $\mathcal{H}_{2}$ approximation error bound.


\section*{Appendix A. Proof of Theorem \ref{theorem-reducedasy}}\label{sec:appendixa}
\setcounter{equation}{0}
\renewcommand\theequation{A.\arabic{equation}}
\begin{pf}
 Firstly, we prove that the reduced second-order model $\bm{\hat{\Sigma}_{rs}}$ as in \eqref{reduced-1} is asymptotically stable with system matrices $K_r$, $D_r$, $F_r$, $H_r$ given in \eqref{reducedasym2}.
It follows from $\hat{P}_{31}>0$, and full column rank of $\hat {P}_{21}$ that $\hat{P}_{21}\hat{P}^{-1}_{31}$ is a full column rank matrix. Since $K, D>0$, we obtain that
$
 D_{r} =\hat{P}^{-1}_{31}\hat{P}^\top_{21}D
\hat{P}_{21}\hat{P}^{-1}_{31}>0$,
 and $
 K_{r} =\hat{P}^{-1}_{31}\hat{P}^\top_{21}K
\hat{P}_{21}\hat{P}^{-1}_{31}>0.
$
According to \cite{Semistable},
the reduced second-order system $\bm{\Sigma_{r}}$ with system matrices given in \eqref{reducedasym2} is asymptotically stable. Note that if there exist matrices $\hat{P}>0$, $K_{r}$, $D_{r}$, $F_{r}$, $H_{r}$ satisfy the following optimization problem:
\begin{subequations}\label{A-1}
\begin{align}\label{A-1a}
  \min_{\hat{P}>0, \gamma>0} \quad& \gamma\\
\label{A-1b}
  s.t.\quad& \tr\left(C_{e}\hat{P}C^\top_{e}\right)<\gamma,   \\
\label{A-1c}
&\hat{P}A^\top_{e}+A_{e}\hat{P}+B_{e}B^\top_{e}<0,\\
\label{A-1d}
  &\hat{P}=\begin{bmatrix}
  \begin{smallmatrix}
 \hat{P}_{11}& \hat{P}_{12}&\hat{P}_{21}& \hat{P}_{22}\\
 \star& \hat{P}_{13}&\hat{P}_{23}&\hat{P}_{24}\\
 \star& \star&\hat{P}_{31}& \hat{P}_{32}\\
 \star& \star&\star&\hat{P}_{33}\\
 \end{smallmatrix}\end{bmatrix}>0,
\end{align}
\end{subequations}
with $A_{e}$, $B_{e}$, $C_{e}$ given in \eqref{AeBe}.
Then, it follows that the approximation error between the original interconnected second-order system $\bm{\Sigma}$ in \eqref{origi-1} and the reduced second-order model $\bm{\Sigma_{r}}$ in \eqref{reduced-1} satisfies the upper bound given in \eqref{errorbound}.
Now, we prove that if there exist matrices $\hat{P}_{11} >0$,  $\hat{P}_{12}$, $\hat{P}_{13} >0$, $\hat{P}_{31} >0$, $\hat{P}_{31} $, $\hat{P}_{21} $, and $X \geq0$, such that the optimization problem \eqref{optimization2} is solvable, then the optimization problem
\eqref{A-1} is also solvable. That is, the solution of optimization problem \eqref{optimization2} is also a solution of the problem \eqref{A-1}.

In the sequel, we prove that the inequalities \eqref{trace1}-\eqref{linearinequ3} are the necessary and sufficient conditions for the  problem \eqref{A-1} when $\hat{P}$ has the form of \eqref{hatP1}. Note that the inequality \eqref{A-1c} can be rewritten as
\begin{equation}\label{A-5}
  G+\sym(K_{1}Y K^\top_{2})<0,
\end{equation}
where
\begin{equation*}
\begin{split}
G&=\begin{bmatrix}\begin{smallmatrix}
\sym(\hat{P}_{12})& -\hat{P}_{11}K-\hat{P}_{12}D+\hat{P}_{13}
&~0&0&0\\
\star& \sym(-K\hat{P}_{12}-D\hat{P}_{13})& ~-K\hat{P}_{21}&~ 0& F\\
\star&\star& 0&0&0\\
\star&\star&\star&0&0\\
\star&\star&\star&\star& -I
\end{smallmatrix}
\end{bmatrix},\\
K_{1}&=\begin{bmatrix}
\begin{smallmatrix}
  0&0\\
  0&0\\
  I&0\\
  0&I\\
  0&0
\end{smallmatrix}
\end{bmatrix},
K_2=\begin{bmatrix}
\begin{smallmatrix}
  \hat{P}_{21}& 0&0\\
  0&0&0\\
  \hat{P}_{31}& -\hat{P}_{31}&0\\
  -\hat{P}_{31}& 2\hat{P}_{31}&0\\
    0& 0&I\\
\end{smallmatrix}
\end{bmatrix},Y=\begin{bmatrix}
\begin{smallmatrix}
 0& I_{r}& 0\\
 -K_{r}& -D_{r}& F_{r}
\end{smallmatrix}
\end{bmatrix}
\end{split}
\end{equation*}
and the orthogonal complements of the matrices $K_1$, $K_2$
are given by
\begin{equation*}
    K^{\perp}_{1} =\begin{bmatrix}
    \begin{smallmatrix}
I& 0&0& 0& 0\\
0& I& 0& 0&0\\
0& 0& 0& 0& I
\end{smallmatrix}
\end{bmatrix},   
K^{\perp}_{2} =\begin{bmatrix}
\begin{smallmatrix}
 I& 0&~-2\hat{P}_{21}\hat{P}^{-1}_{31}& ~-\hat{P}_{21}\hat{P}^{-1}_{31}& 0\\
0&I &0& 0& 0
\end{smallmatrix}
\end{bmatrix}.
\end{equation*}
According to the Finsler's lemma, the inequality \eqref{A-5} is equivalent to
  $K^{\perp}_{1}G(K^{\perp}_{1})^\top<0$, $K^{\perp}_{2}G(K^{\perp}_{2})^\top<0$,
where the first inequality is equivalent to
\begin{equation*}
   \begin{bmatrix}
    \sym(\hat{P}_{12})&~ \Pi_{12}\\
 \star& \Pi_{22}
   \end{bmatrix}<0,
\end{equation*}
as given in \eqref{linearinequ1}, and the second inequality is equivalent to $\Phi<0$, as given in \eqref{linearinequ2}.

Next, we prove that \eqref{trace1} is a necessary condition of inequality \eqref{A-1b}. Suppose that $R=R^\top>0$ satisfies $R-C_{e}\hat{P}C^\top_{e}>0$. Therefore, $\tr(R)<\gamma^{2}$ implies $\tr(C_{e}\hat{P}C^\top_{e})<\gamma^{2}$. By using Schur complement, $R-C_{e}\hat{P}C^\top_{e}>0$ is equivalent to
  $\begin{bmatrix}
  \begin{smallmatrix}
   R& ~C_{e}\hat{P}\\
   \star& ~\hat{P}
   \end{smallmatrix}
  \end{bmatrix}>0$,
which can be rewritten as
\begin{equation}\label{A-7}
  \Omega-\sym(\Upsilon\begin{bmatrix}
   H_{r}&~0
  \end{bmatrix}Z^\top)>0,
\end{equation}
where
\begin{equation*}
\begin{split}
 \Omega&=\begin{bmatrix}
 \begin{smallmatrix}
 R&~ H\hat{P}_{11}&~ H\hat{P}_{12}&~
H\hat{P}_{21}& 0\\
 \star& \hat{P}_{11}& \hat{P}_{12}& \hat{P}_{21}& 0\\
 \star&\star& \hat{P}_{13}& 0& 0\\
 \star&\star&\star&  \hat{P}_{31}& -\hat{P}_{31}\\
 \star&\star&\star&\star& 2\hat{P}_{31}
 \end{smallmatrix}
\end{bmatrix},~
\Upsilon=\begin{bmatrix}
\begin{smallmatrix}
  I\\
  0\\
  0\\
  0\\
  0
\end{smallmatrix}
\end{bmatrix}, \\ Z&=\begin{bmatrix}
\begin{smallmatrix}
  0&\hat{P}^\top_{21}& 0&~\hat{P}_{31}& ~-\hat{P}_{31}\\
 0& 0&0&~ -\hat{P}_{31}& ~2\hat{P}_{31}
\end{smallmatrix}
\end{bmatrix}^\top.
\end{split}
\end{equation*}
The orthogonal complements of matrices
$\Upsilon$, $Z$ are
\begin{equation*}
\Upsilon^{\perp}=\begin{bmatrix}
\begin{smallmatrix}
0& I&0&0&0\\
0&0& 0& 2I& I
\end{smallmatrix}
\end{bmatrix},  \
Z^{\perp}=\begin{bmatrix}
\begin{smallmatrix}
I&0&0&0&0\\
0& I& 0& ~ -2\hat{P}_{21}\hat{P}^{-1}_{31}& ~-\hat{P}_{21}\hat{P}^{-1}_{31}\\
0& 0&I& 0&0
\end{smallmatrix}
\end{bmatrix}.
\end{equation*}
According to the Finsler's lemma, \eqref{A-7} is equivalent to
$\Upsilon^{\perp}\Omega (\Upsilon^{\perp})^\top>0$, $Z^{\perp}\Omega (Z^{\perp})^\top>0$,
which can be rewritten as inequality \eqref{linearinequ3} and
\begin{equation*}
\begin{bmatrix}
\begin{smallmatrix}
  R&~ H\hat{P}_{11}-2HX& ~ H\hat{P}_{12}\\
  \star& \hat{P}_{11}-2X&~ \hat{P}_{12}\\
  \star&\star&\hat{P}_{13}
\end{smallmatrix}
\end{bmatrix}>0.
\end{equation*}
The above inequality leads to
\begin{equation}\label{A-8}
  \begin{bmatrix}
  \begin{smallmatrix}
   R&~ ~H\hat{P}_{11}-2HX\\
   \star& \hat{P}_{11}-2X
   \end{smallmatrix}
  \end{bmatrix}>0.
\end{equation}
By using Schur complement, \eqref{A-8} is equivalent to
\begin{equation*}
R-H(\hat{P}_{11}-2X)
  H^\top>0.
\end{equation*}
Thus, $\tr\left(H
(\hat{P}_{11}-2X)
H^\top\right)<\gamma^{2}$ appears as a necessary condition to satisfy
$
 H(\hat{P}_{11}-2X)
 H^\top<R$,  and $\tr(R)<\gamma^{2}.
$
Note that the rank of $\hat{P}_{21}\in\mathbb{R}^{n\times r}$ could not exceed $r$ since the projection matrix $\hat{P}_{21}\hat{P}^{-1}_{31}$ must have full column rank, that is, $\rank(\hat{P}_{21}\hat{P}^{-1}_{31})=r$. Therefore, the rank of $X=\hat{P}_{21}\hat{P}^{-1}_{31}\hat{P}^\top_{21}$ satisfies $\rank(X)\leq r$.
This completes the proof of Theorem \ref{theorem-reducedasy}.
	\hfill $\openbox$
\end{pf}

\bibliographystyle{IEEEtran}
\bibliography{second-order}
\end{document}